\documentclass[12pt,a4paper]{article}
\usepackage[utf8]{inputenc}
\usepackage{amsmath}
\usepackage{amsfonts}
\usepackage{amssymb}
\usepackage{amsthm}
\usepackage{authblk}
\usepackage{inputenc}
\allowdisplaybreaks

\begin{document}
\title{Properties of counterexample to Robin hypothesis}
\author{Xiaolong Wu}
\affil{Ex. Institute of Mathematics, Chinese Academy of Sciences}
\affil{xwu622@comcast.net}
\date{January 7, 2019}
\maketitle

\begin{abstract}

    Let $G(n)=\sigma (n)/(n \log \log n )$. Robin made hypothesis that $G(n)<e^\gamma$ for all integer $n>5040$. If there exists counterexample to Robin hypothesis, then there must exist finite number of counterexamples $n>5040$ such that $G(n)$ attains largest value. This article studies various properties of such number.
\end{abstract}
\begin{center}\textbf{ \large Introduction}
\end{center}

Robin made a hypothesis [Robin 1984] that the Robin's inequality 
\begin{equation}\tag{RI}
\sigma (n)<e^\gamma n \log \log n, 
\end{equation}
 holds for all integers $n>5040$. Here $\sigma (n)=\sum_{d|n} d$ is the divisor sum function, $\gamma $ is the Euler-Mascheroni constant, log is the nature logarithm.

 For calculation convenience, we define \[\rho (n):=\frac{\sigma (n)}{n}.\]
Then Robin’s inequality can also be written as 
\begin{equation}\tag{RI}
\rho (n)<e^\gamma \log \log n .
\end{equation}

Define
\begin{equation}\notag
G(n):=\frac{\rho (n)}{\log \log n} .
\end{equation}
Then Robin’s inequality can also be written as 
\begin{equation}\tag{RI}
G(n)<e^\gamma.
\end{equation}

    Let $N>5040$ be an integer. Write the factorization of N as 
\[N=\prod_{i=1}^{r} p_i^{a_i},\]
where $p_i$ are in increasing orders, $p_r$ is the largest prime factor of N.

 According to $[Morrill; Platt\, 2018]$, (RI) holds for all integers  $n, 5040<n\leq 10^{(10^{13})}$. So, we assume $N>10^{(10^{13})}$.
 
    By Grönwall’s theorem, {[}Broughan 2017{]} Theorem 9.2, if there exist counterexamples of Robin hypothesis, then there must exist finite number of counterexamples $n>5040$ such that G(n) attains largest value. We call such an n a largest G-value (abbreviate LG) number. 
    
    This article proves the following properties of LG numbers. Assume N is an LG number. Then\\
	1) N is colossally abundant.\\
	2) $p_r<log N$.\\
	3) $p_r$ is the largest prime below N.\\
	4)\[\, a_i\leq \left\lfloor \frac{\log (kp_r )}{\log p_i}\right\rfloor,\,when\, ((k+1) p_r )^{1/((k+1)}<p_i\leq (kp_r )^{1/k}, \forall\, k\geq 1.\]
	5)\[\, a_i\geq \left\lfloor \frac{\log p_r}{\log p_i} \right\rfloor \forall \, i\leq r.\]
	6)\[\, \log N>p_r+\frac{1}{2}\log p_r+\frac{1}{2}-\frac{1}{2\log p_r}.\]
	7) Let p be the smallest prime above $\log N$, then \[\log N<p-\frac{1}{2}\log p+\frac{1}{2}-\frac{1}{2}\log p+\frac{1}{(\log p)(\log p+1)}.\]
	8)\[\, G(N)<e^\gamma+\frac{0.00995}{(\log \log N)^2}.  \]
	9)\[\, p_r>\log N\left(1-\frac{0.005587}{\log \log N}\right)\, and\, \log N\leq p_r\left(1+\frac{0.005589}{\log p_r}\right).\]

{\bfseries Version Notes:}\\
2019-02-13 version 2. Added two theorems. They are reverse of theorems 6 and 7.\\
Theorem 10. $G(N)>G(N/p)$ if 
\[\log N>p+\frac{\log p}{2}+\frac{1}{2}-\frac{1}{2\log p}+\frac{1}{(\log p)(\log p+1}.\]
Theorem 11. $G(N)>G(Np)$ if 
\[\log N<p-\frac{1}{2}\log p+\frac{1}{2}-\frac{1}{\log p+1}\]
\begin{center}
\textbf{ \large Main Content}
\end{center}

\noindent {\bfseries Theorem 1.}
\textit{Let N be an LG number, then N is colossally abundant.
}
\begin{proof}
By Proposition 1 of {[}Robin 1984{]}, N is between two adjacent colossally numbers  $n_i$ and $n_{i+1}$ for some integer i. We have 
\begin{equation}\notag
G(N)\leq \max(G(n_i ),G(n_{i+1}).
\end{equation}
By maximality of $G(N)$, the equal sign must hold. By strict convexity of $x\rightarrow \epsilon x-\log \log x\, (x>1)$, we must have $N=n_i$ or $N=n_{i+1}$.
\end{proof}

\noindent {\bfseries Theorem 2.}
\textit{Let N be an LG number. Then $p_r<\log N$.
}
\begin{proof} Write $p:=p_r$. By Theorem 1, we know N is colossally abundant, so the exponent of p in N is 1. We have
\begin{align*}
\frac{G(N)}{G(N/p)}&=\frac{\rho(N)\log \log (N/p)}{\rho(N/p)\log \log N}\\
&=\frac{\log(\log N-\log p)}{\log \log N}\left(1+\frac{1}{p}\right)\\
&=\frac{\log \log N+\log \left(1-\frac{\log p}{\log N}\right)}{\log \log N}\left(1+\frac{1}{p}\right)\\
&<\left(1-\frac{\log p}{\log N\log \log N}\right)\left(1+\frac{1}{p}\right).\tag{2.1}
\end{align*}
If $p\geq \log N$, we would have
\begin{equation}\tag{2.2}
\frac{G(N)}{G(N/p)}<1+\frac{log N\log \log N-p \log p-\log p}{p\log N\log \log N}<1.
\end{equation}
That is, $G(N)<G(N/p)$, which contradicts to the maximality of N.
\end{proof}

\noindent {\bfseries Theorem 3.}
\textit{Let N be an LG number. Then $p_r$ must be the largest prime below $\log N$.
}
\begin{proof}
We know $p_r<\log N$ by Theorem 2. Assume there exists a prime p such that $p_r<p<\log N$. We will derive a contradiction.  Compare $G(N)$ and $G(Np)$, we have 
\begin{align*}
\frac{G(N)}{G(Np)}&=\frac{\rho(N)\log \log (Np)}{\rho(Np)\log \log N}\\
&=\frac{\log(\log N+\log p)}{\log \log N}\left(\frac{p}{p+1}\right)\\
&=\frac{\log \log N+\log \left(1+\frac{\log p}{\log N}\right)}{\log \log N}\left(\frac{p}{p+1}\right)\\
&<\left(1+\frac{\log p}{\log N\log \log N}\right)\left(\frac{p}{p+1}\right).\tag{3.1}
\end{align*}
Since $p<\log N$, we have
\begin{equation}\tag{3.2}
\frac{G(N)}{G(Np)}<\left(1+\frac{\log p}{p\log p}\right)\left(\frac{p}{p+1}\right)=1.
\end{equation}
That means $G(N)<G(Np)$, which contradicts to the maximality of N.
\end{proof}

Recall the construction of a colossally abundant number $N_\epsilon$ from a given parameter $\epsilon>0$, cf. {[}EN 1975{]} Proposition 4 or {[}Broughan 2017{]} Section 6.3. Define
\begin{equation}\notag
N_\epsilon :=\prod_p p^{a_p(\epsilon)},\quad a_p(\epsilon):=\left\lfloor \frac{\log ((p^{1+\epsilon}-1)/(p^\epsilon-1))}{\log p}\right\rfloor-1.
\end{equation}
Let $k\geq 1$ be an integer, $x_k$ be the solution of
\begin{equation}\notag
F(x,k):=\frac{\log (1+1/(x+x^2+\cdots+x^k))}{\log x}=\epsilon.
\end{equation}
Then one can show that
\begin{equation}\notag
a_p(\epsilon)=\begin{cases}
k, & \text{if $x_{k+1}<p\leq x_k,\,k\geq 1$}\\
0, & \text{if $p>x_1$.}
\end{cases}
\end{equation}

\noindent {\bfseries Theorem 4.}
\textit{Let $\epsilon>0$ be a parameter, $N_\epsilon$ be the colossally number constructed from $\epsilon$,  $p\geq 3299$ be the largest prime factor of $N_\epsilon$. Then 
\begin{equation}\tag{4.1}
x_k<(kp)^{1/k}, \quad \forall\, k \geq 2.
\end{equation}
}
\begin{proof} This is an improvement based on Lemma 1 of {[}CNS 2012{]}, which proved $x_k<(kx_1)^{1/k},\forall\, k\geq2$. Since the function $t\rightarrow F(t,k)$ is strictly decreasing on $1<t<\infty$, to prove that $x_k<z:=(kp)^{1/k}$, it suffices to show $F(z,k)<F(x_k,k)$. Since $F(x_k,k)=\epsilon =F(x_1,1)$, this reduces to showing $F(z,k) < F(x_1,1)$.
\begin{align*}
F(z,k)&=\log \left(1+\frac{1}{z+z^2+\cdots+z^k}\right)\frac{1}{\log z}\\
&<\frac{1}{(z+z^2+\cdots+z^k)\log z}<\frac{k}{(k-1+z^k)\log kp}\\
&\leq \frac{k}{\left(\frac{k}{2}+z^k\right)\log kp}=\frac{1}{\left(p+\frac{1}{2}\right)\log kp}\\
&<\log \left(1+\frac{1}{p}\right)\frac{1}{\log kp}.\tag{4.2}
\end{align*}
We need to show
\begin{equation}\tag{4.3}
\log \left(1+\frac{1}{p}\right)\frac{1}{\log kp}<F(x_1,1)=\log \left(1+\frac{1}{x_1}\right)\frac{1}{\log x_1},
\end{equation}
that is
\begin{equation}\tag{4.4}
\frac{\log \left(1+\frac{1}{p}\right)}{\log \left(1+\frac{1}{x_1}\right)}<\frac{\log kp}{\log x_1}.
\end{equation}

Write
\begin{equation}\notag
g(t):=t \log \left(1+\frac{1}{t}\right).
\end{equation}
Take derivative
\begin{align*}
g'(t)&=\log \left(1+\frac{1}{t}\right)+\frac{t}{1+\frac{1}{t}}\cdot \frac{-1}{t^2}=\sum_{j=1}^\infty \frac{(-1)^{j-1}}{jt^j}-\frac{1}{t+1}\\
&>\frac{1}{t}-\frac{1}{2t^2}-\frac{1}{t+1}=\frac{2t^2+2t-t-1-2t^2}{2t^2(t+1)}=\frac{t-1}{2t^2(t+1)}>0,\tag{4.5}
\end{align*}
for $t>1$. Hence $g(t)$ strictly increases, and
\begin{equation}\tag{4.6}
\frac{\log \left(1+\frac{1}{p}\right)}{\log \left(1+\frac{1}{x_1}\right)}<\frac{x_1}{p}.
\end{equation}
So in view of (4.4), it suffices to prove
\begin{equation}\tag{4.7}
\frac{x_1}{p}<\frac{\log kp}{\log x_1}.
\end{equation}
By Proposition 5 of {[}Dusart 1998{]}, for all $j\geq 463, (p_463=3299)$, we have
\begin{equation}\tag{4.8}
p_{j+1}\leq p_j\left(1+\frac{1}{2(\log p_j)^2}\right).
\end{equation}
Theorem assumes $p\geq3299$. Since p is the largest prime $\leq x_1$, we must have
\begin{equation}\tag{4.9}
x_1< p\left(1+\frac{1}{2(\log p)^2}\right).
\end{equation}
(4.7) becomes
\begin{align*}
\frac{p \log kp}{x_1 \log x_1}&>\frac{p(\log p+\log k)}{p\left(1+\frac{1}{2(\log p)^2}\right)\log \left(p\left(1+\frac{1}{2(\log p)^2}\right)\right)}\\
&=\frac{\log p+\log k}{\left(1+\frac{1}{2(\log p)^2}\right)\left(\log p+\log \left(1+\frac{1}{2(\log p)^2}\right)\right)}\\
&>\frac{\log p+\log k}{\left(1+\frac{1}{2(\log p)^2}\right)\left(\log p+\frac{1}{2(\log p)^2}\right)}\\
&=\frac{\log p+\log k}{\log p+\frac{1}{2\log p}+\frac{1}{2(\log p)^2}+\frac{1}{4(\log p)^4}}.\tag{4.10}
\end{align*}
Since
\begin{equation}\tag{4.11}
\frac{1}{2\log p}+\frac{1}{2(\log p)^2}+\frac{1}{4(\log p)^4}<\log k\quad \forall \, p\geq 5,\,k\geq 2,
\end{equation}
we have
\begin{equation}\tag{4.12}
\frac{p \log kp}{x_1 \log x_1}>1,
\end{equation}
i.e. (4.7) holds.
\end{proof}

{\bfseries Definition 1.} Now we construct a lower bound curve L for the exponents. Define
\begin{equation}\tag{D1.1}
L(p_i )=L_{p_r}(p_i ):=\left\lfloor \frac{\log p_r}{\log p_i}\right\rfloor\quad  for\, i\leq r.
\end{equation}

\noindent {\bfseries Theorem 5.}
\textit{ Let $N>10^{\left(10^{13}\right)}$ be an LG number. Then $a_i \geq L(p_i )$. 
}
\begin{proof}
As N being a colossally abundant number, we know $a_r=1=L(p_r)$. Assume $a_s< L(p_i )$ for some index $s<r$. We will derive a contradiction. Define
\[N_1:=(p_s/p_r)N.\]
Then $\log N-\log N_1=\log p_r-\log p_s$. $p_s<p_r$ means $N_1<N$. $a_s<L(p_s )=\left\lfloor\frac{\log p_r}{\log p_s} \right\rfloor$ means $a_s+1\leq \left\lfloor \frac{\log p_r}{\log p_s} \right\rfloor \leq \frac{\log p_r}{\log p_s}$. Hence $p_s^{a_s+1}\leq p_r$ and 
\begin{equation}\tag{5.1}
\log p_s \leq \frac{1}{a_s+1} \log p_r.
\end{equation}
It is easy to deduce
\begin{align*}
\frac{G(N)}{G(N_1)}&=\frac{\rho(N)\log \log (N_1)}{\rho(N_1)\log \log N}\\
&=\frac{\log(\log N-\log p_r+\log p_s)}{\log \log N}\left(\frac{p_s-p_s^{-a_s}}{p_s-p_s^{-a_s-1}}\right)\left(\frac{p_r+1}{p_r}\right)\\
&\leq \left(1-\frac{\log p_r-\frac{1}{a_s+1}\log p_r}{\log N \log \log N} \right)\left(\frac{p_s-p_s^{-a_s}}{p_s-p_s^{-a_s-1}}\right)\left(1+\frac{1}{p_r}\right)\\
&=\left(1-\left(\frac{a_s}{a_s+1}\right)\frac{\log p_r}{\log N\log \log N}\right)\left(\frac{p_s-p_s^{-a_s}}{p_s-p_s^{-a_s-1}}\right)\left(1+\frac{1}{p_r}\right).\tag{5.2}
\end{align*}
\begin{equation}\tag{5.3}
\frac{p_s-p_s^{-a_s}}{p_s-p_s^{-a_s-1}}=1-\frac{1}{p_s^{a_s+1}+p_s^{a_s}+\cdots+1}.
\end{equation}
By Proposition 5 of {[}Dusart 1998{]}, for all $j\geq 463$, $(p_{463}=3299)$, we have
\[ p_{j+1}\leq p_j \left(1+\frac{1}{2(\log p_j)^2}\right).\]
By Theorem 3, $p_r$ is the largest prime below $\log N$, so
\begin{equation}\tag{5.4}
p_r>\log N \left(1-\frac{1}{2(\log p_r)^2}\right).
\end{equation}
We have, noting $N>10^{\left(10^{13}\right)}$,
\begin{equation}\tag{5.5}
\log N<cp_r, \quad c:= \left(1-\frac{1}{2(\log (2.3\times 10^{13}))^2}\right)=1.000528\cdots.
\end{equation}
Since $\log (cp_r)<c \log p_r$, (5.2) can be simplified to
\begin{align*}
\frac{G(N)}{G(N_1)}&<\left(1-\left(\frac{a_s}{a_s+1}\right)\frac{\log p_r}{(cp_r)\log (cp_r)}\right)\left(\frac{p_s-p_s^{-a_s}}{p_s-p_s^{-a_s-1}}\right)\left(1+\frac{1}{p_r}\right)\\
&<\left(1-\left(\frac{a_s}{a_s+1}\right)\frac{1}{c^2p_r}\right)\left(\frac{p_s-p_s^{-a_s}}{p_s-p_s^{-a_s-1}}\right)\left(1+\frac{1}{p_r}\right).\tag{5.6}
\end{align*}
Now we split the proof into two cases.\\
\textbf{Case 1)} $a_s=1$. We have in this case
\begin{equation}\tag{5.7}
1-\left(\frac{a_s}{a_s+1}\right)\frac{1}{c^2p_r}<1-\frac{1}{2c^2p_r}<1-\frac{0.49}{p_r}
\end{equation}
\begin{equation}\tag{5.8}
p_s^{a_s+1}+p_s^{a_s}+\cdots+1=p_s^2+p_s+1\leq\frac{7}{4}p_s^2,
\end{equation}
\begin{equation}\tag{5.9}
\frac{p_s-p_s^{-a_s}}{p_s-p_s^{-a_s-1}}=1-\frac{1}{p_s^2+p_s+1}\leq 1-\frac{4}{7p_s^2}<1-\frac{0.57}{p_r}.
\end{equation}
Substitute (5.7) and (5.9) in to (5.6), we get 
\begin{equation}\tag{5.10}
\frac{G(N)}{G(N_1)}<\left(1-\frac{0.49}{p_r}\right)\left(1-\frac{0.57}{p_r}\right)\left(1+\frac{1}{p_r}\right)<1,
\end{equation}
which contradicts to the maximality of N.\\
\textbf{Case 2)} $a_s>1$. We have
\begin{equation}\tag{5.11}
1-\left(\frac{a_s}{a_s+1}\right)\frac{1}{c^2p_r}<1-\frac{2}{3c^2p_r}<1-\frac{0.66}{p_r}.
\end{equation}
\begin{equation}\tag{5.12}
\frac{p_s-p_s^{-a_s}}{p_s-p_s^{-a_s-1}}=1-\frac{1}{p_s^{a_+1}+p_s^{a_s}+\cdots+1}<1-\frac{1}{2p_s^{a_s+1}}<1-\frac{0.50}{p_r}.
\end{equation}
Substitute (5.11) and (5.12) in to (5.6), we get
\begin{equation}\tag{5.13}
\frac{G(N)}{G(N_1)}<\left(1-\frac{0.66}{p_r}\right)\left(1-\frac{0.50}{p_r}\right)\left(1+\frac{1}{p_r}\right)<1,
\end{equation}
which contradicts to the maximality of N.
\end{proof}

\noindent {\bfseries Lemma 1.}
\textit{ Let N be an integer, p be a prime factor of N with exponent 1. Write $\log N=p+\frac{1}{2}\log p+d$.
Then $G(N)>G(N/p)$ if and only if
\begin{align*}
&\frac{1}{2}(\log p)^2+d \log p+d+\sum_{k=1}^\infty \frac{(-1)^{k-1}\left(\frac{1}{2}\log p+d\right)^{k+1}}{k(k+1)p^k}\\
>\frac{1}{2} \log p&+\frac{p(\log p)^2}{2\log N}+\frac{(\log p)^2}{2\log N}+(p+1)\sum_{k=1}^\infty \frac{(\log p)^{k+2}}{(k+2)(\log N)^{k+1}}.\tag{L1.1}
\end{align*}
}
\begin{proof}
Substitute $\log N$
\begin{align*}
\log N &\log \log N=\left(p+\frac{1}{2}\log p+d\right)\log \left(p+\frac{1}{2}\log p+d\right)\\
&=\left(p+\frac{1}{2}\log p+d\right)\left(\log p+\log \left(1+\frac{\frac{1}{2}\log p+d}{p}\right)\right)\\
&=p\log p+\frac{1}{2}(\log p)^2+d \log p\\
&\quad+\left(p+\frac{1}{2}\log p+d\right)\log \left(1+\frac{\frac{1}{2}\log p+d}{p}\right)\\
&=p\log p+\frac{1}{2}(\log p)^2+d \log p\\
&\quad+\left(p+\frac{1}{2}\log p+d\right)\sum_{k=1}^\infty \frac{(-1)^{k-1}\left(\frac{1}{2}\log p+d\right)^k}{kp^k}\\
&=p\log p+\frac{1}{2}(\log p)^2+d \log p\\
&\quad+\sum_{k=1}^\infty \frac{(-1)^{k-1}\left(\frac{1}{2}\log p+d\right)^k}{kp^{k-1}}+\sum_{k=1}^\infty \frac{(-1)^{k-1}\left(\frac{1}{2}\log p+d\right)^{k+1}}{kp^k}\\
&=p\log p+\frac{1}{2}(\log p)^2+d \log p+\frac{1}{2}\log p+d\\
&\quad+\sum_{k=1}^\infty \left(\frac{(-1)^k\left(\frac{1}{2}\log p+d\right)^{k+1}}{(k+1)p^k}+\frac{(-1)^{k-1}\left(\frac{1}{2}\log p+d\right)^{k+1}}{kp^k}\right)\\
&=p\log p+\frac{1}{2}(\log p)^2+d \log p+\frac{1}{2}\log p+d\\
&\quad+\sum_{k=1}^\infty \frac{(-1)^{k-1}\left(\frac{1}{2}\log p+d\right)^{k+1}}{k(k+1)p^k}\tag{L1.2}
\end{align*}
Compare $G(N)$ and $G(N/p)$, we have
\begin{align*}
\frac{G(N)}{G(N/p)} &=\frac{\rho(N)\log \log (N/p)}{\rho(N/p)\log \log N}\\
&=\frac{\log (\log N-\log p)}{\log \log N}\left(1+\frac{1}{p}\right)\\
&=\left(1+\frac{\log \left(1-\frac{\log p}{\log N}\right)}{\log \log N}\right)\left(1+\frac{1}{p}\right).\tag{L1.3}
\end{align*}
Therefore,
\begin{align*}
&\quad \quad G(N)>G(N/p)\\
&\Longleftrightarrow \left(1+\frac{\log \left(1-\frac{\log p}{\log N}\right)}{\log \log N}\right)\left(1+\frac{1}{p}\right)>1\\
&\Longleftrightarrow 1+\frac{\log \left(1-\frac{\log p}{\log N}\right)}{\log \log N}>\left(1+\frac{1}{p}\right)^{-1}=1-\frac{1}{p+1} \\
&\Longleftrightarrow -\frac{1}{\log \log N}\sum_{k=1}^\infty \frac{1}{k}\left(\frac{\log p}{\log N}\right)^k>-\frac{1}{p+1} \\
&\Longleftrightarrow \frac{\log p}{\log N \log \log N}\left(1+\sum_{k=2}^\infty \frac{1}{k}\left(\frac{\log p}{\log N}\right)^{k-1}\right)<\frac{1}{p+1} \\
&\Longleftrightarrow (p+1)\log p\left(1+\sum_{k=1}^\infty \frac{1}{k+1}\left(\frac{\log p}{\log N}\right)^k\right)<\log N \log \log N \tag{L1.4}
\end{align*}
Compare (L1.2) and (L1.4), we see that $G(N)>G(N/p)$ if and only if
\begin{align*}
p\log p&+\frac{1}{2}(\log p)^2+d \log p+\frac{1}{2}\log p+d\\
&\quad +\sum_{k=1}^\infty \frac{(-1)^{k-1}\left(\frac{1}{2}\log p+d\right)^{k+1}}{k(k+1)p^k}\\
&>(p+1)\log p\left(1+\frac{\log p}{2\log N}+\sum_{k=2}^\infty \frac{1}{k+1}\left(\frac{\log p}{\log N}\right)^k\right)\\
&=p\log p+\log p+\frac{p(\log p)^2}{2\log N}+\frac{(\log p)^2}{2\log N}\\
&\quad +(p+1)\sum_{k=1}^\infty \frac{(\log p)^{k+2}}{(k+2)(\log N)^{k+1}}.\tag{L1.5}
\end{align*}
\end{proof}

\noindent {\bfseries Theorem 6. }
\textit{ Let $N>10^{(10^{13})}$ be an LG number. Then 
\begin{equation}\tag{6.1}
\log N>p_r+\frac{1}{2}\log p_r+\frac{1}{2}-\frac{1}{2\log p_r}.
\end{equation}
}
\begin{proof} 
Write $p:=p_r$. By Theorem 3, p is the largest prime below $\log N$. Write $\log N=p+\frac{1}{2}\log p+d$, where d is a to-be-determined expression. By Lemma 1, $G(N)>G(N/p)$ if and only if
\begin{align*}
&\frac{1}{2}(\log p)^2+d \log p+d+\sum_{k=1}^\infty \frac{(-1)^{k-1}\left(\frac{1}{2}\log p+d\right)^{k+1}}{k(k+1)p^k}\\
>\frac{1}{2} \log p&+\frac{p(\log p)^2}{2\log N}+\frac{(\log p)^2}{2\log N}+(p+1)\sum_{k=1}^\infty \frac{(\log p)^{k+2}}{(k+2)(\log N)^{k+1}}\tag{6.2}
\end{align*}
This implies
\begin{equation}\tag{6.3}
\frac{1}{2}(\log p)^2+d \log p+d+\frac{\left(\frac{1}{2}\log p+d)\right)^2}{2p}>\frac{\log p}{2}+\frac{p(\log p)^2}{2\log N}. 
\end{equation}
Since p is the largest prime below $\log N$, by Proposition 5.4 of {[}Dusart 2018{]}, for $p\geq 89\,693$ we have
\begin{equation}\tag{6.4}
p>\log N\left(1-\frac{1}{(\log p)^3}\right),
\end{equation}
\begin{equation}\tag{6.5}
\frac{p}{2\log N}>\frac{1}{2}\left(1-\frac{1}{(\log p)^3}\right).
\end{equation}
Since $N>10^{(10^{13})}$, $\log N>(\log 10 )\times 10^{13}$, the last term on left of (6.3) is in order of $10^{-13}(\log p)^2$ and can be absorbed by rounding: the numerator 1 in (6.4) was rounded from 0.998. We can concentrate on main terms. $G(N)>G(N/p)$ implies
\begin{equation}\tag{6.6}
\frac{1}{2}(\log p)^2+d\log p+d>\frac{\log p}{2}+\frac{(\log p)^2}{2}\left(1-\frac{1}{(\log p)^3}\right).
\end{equation}
Hence
\[d(\log p+1)>\frac{1}{2}\log p-\frac{1}{2\log p}.\]
\begin{equation}\tag{6.7}
d>\frac{\log p-\frac{1}{\log p}}{2(\log p+1)}=\frac{1-\frac{1}{(\log p)^2}}{2\left(1+\frac{1}{\log p}\right)}=\frac{1}{2}\left(1-\frac{1}{\log p}\right).
\end{equation}
\end{proof}

\noindent {\bfseries Lemma 2. }
\textit{ Let $N>5040$ be an integer. $p<N$ be a prime. Assume p does not divide N. Write $\log N=p-\frac{1}{2}\log p+d$. Then $G(N)>G(Np)$ if and only if
\begin{align*}
&p \sum_{k=1}^\infty \frac{(-1)^k(\log p)^{k+1}}{(k+1)(\log N)^k}\\
&>-\frac{1}{2}(\log p)^2+d\log p-\frac{1}{2}\log p+d
+\sum_{k=1}^\infty \frac{\left(\frac{1}{2}\log p-d\right)^{k+1}}{k(k+1)p^k}.\tag{L2.1}
\end{align*}
}
\begin{proof} Substitute $\log N$
\begin{align*}
\log N\log \log N&=\left(p-\frac{1}{2}\log p+d\right)\log \left(p-\frac{1}{2}\log p+d\right)\\
&=\left(p-\frac{1}{2}\log p+d\right)\left(\log p+\log\left(1-\frac{\frac{1}{2}\log p-d}{p}\right)\right)\\
&=p\log p-\frac{1}{2}(\log p)^2+d\log p\\
&\quad -\left(p-\frac{1}{2}\log p+d\right)\sum_{k=1}^\infty \frac{1}{k}\left(\frac{\frac{1}{2}\log p-d}{p}\right)^k\\
&=p\log p-\frac{1}{2}(\log p)^2+d\log p\\
&\quad -\left(\sum_{k=1}^\infty \frac{1}{k}\frac{\left(\frac{1}{2}\log p-d\right)^k}{p^{k-1}}-\sum_{k=1}^\infty \frac{1}{k}\frac{\left(\frac{1}{2}\log p-d\right)^{k+1}}{p^k}\right)\\
&=p\log p-\frac{1}{2}(\log p)^2+d\log p-\frac{1}{2}\log p+d\\
&\quad -\left(\sum_{k=2}^\infty \frac{\left(\frac{1}{2}\log p-d\right)^k}{kp^{k-1}}-\sum_{k=1}^\infty \frac{\left(\frac{1}{2}\log p-d\right)^{k+1}}{kp^k}\right)\\
&=p\log p-\frac{1}{2}(\log p)^2+d\log p-\frac{1}{2}\log p+d\\
&\quad +\sum_{k=1}^\infty \frac{\left(\frac{1}{2}\log p-d\right)^{k+1}}{k(k+1)p^k}\tag{L2.2}
\end{align*}
Compare $G(N)$ and $G(Np)$, we have
\begin{align*}
\frac{G(N)}{G(Np)} &=\frac{\rho(N)\log \log (Np)}{\rho(Np)\log \log N}\\
&=\frac{\log (\log N+\log p)}{\log \log N}\left(\frac{p}{p+1}\right)\\
&=\left(1+\frac{\log \left(1+\frac{\log p}{\log N}\right)}{\log \log N}\right)\left(\frac{p}{p+1}\right).\tag{L2.3}
\end{align*}
Therefore,
\begin{align*}
&G(N)>G(Np)\\
\Longleftrightarrow  &\left(1+\frac{\log \left(1+\frac{\log p}{\log N}\right)}{\log \log N}\right)\left(\frac{p}{p+1}\right)>1\\
\Longleftrightarrow &1+\frac{\log \left(1+\frac{\log p}{\log N}\right)}{\log \log N}>\left(\frac{p}{p+1}\right)^{-1}=1+\frac{1}{p}\\
\Longleftrightarrow &\frac{\log p}{\log N \log \log N}\left(1+\sum_{k=2}^\infty \frac{(-1)^{k-1}}{k}\left(\frac{\log p}{\log N}\right)^{k-1}\right)>\frac{1}{p}\\
\Longleftrightarrow  &p\log p \left(1+\sum_{k=2}^\infty \frac{(-1)^{k-1}}{k}\left(\frac{\log p}{\log N}\right)^{k-1}\right)>\log N \log \log N.\tag{L2.4} 
\end{align*}
Combine (L2.2) and (L2.4), we get
$G(N)>G(Np)$ if and only if 
\begin{align*}
p \log p &\left(1+\sum_{k=1}^\infty \frac{(-1)^k}{k+1}\left(\frac{\log p}{\log N}\right)^k\right)\\
&>p\log p-\frac{1}{2}(\log p)^2+d\log p-\frac{1}{2}\log p+d\\
&\quad +\sum_{k=1}^\infty \frac{\left(\frac{1}{2}\log p-d\right)^{k+1}}{k(k+1)p^k}.\tag{L2.5}
\end{align*}
\end{proof}

\noindent {\bfseries Theorem 7. }
\textit{ Let  $N>5040$ be an integer, p be the prime just above $\log N$. Assume $G(N)>G(Np)$ Then
\begin{equation}\tag{7.1}
\log N<p-\frac{1}{2}\log p+\frac{1}{2}-\frac{1}{2\log p}+\frac{1}{\log p(\log p+1)}.
\end{equation}}
\begin{proof}
Write $\log N=p-\frac{1}{2}\log p+d$, where d is a to-be-determined expression. By Lemma 2, $G(N)>G(Np)$ if and only if 
\begin{align*}
&p \sum_{k=1}^\infty \frac{(-1)^k(\log p)^{k+1}}{(k+1)(\log N)^k}\\\
&>-\frac{1}{2}(\log p)^2+d\log p-\frac{1}{2}\log p+d+\sum_{k=1}^\infty \frac{\left(\frac{1}{2}\log p-d\right)^{k+1}}{k(k+1)p^k}.\tag{7.2}
\end{align*}
So the theorem assumption $G(N)>G(Np)$ implies
\begin{equation}\tag{7.3}
\frac{1}{2}(\log p)^2+\frac{1}{2}\log p-\frac{p(\log p)^2}{2\log N}+\frac{p(\log p)^3}{3(\log N)^2}>d\log p+d.
\end{equation}
Since $p>\log N$, we can replace $\log N$ with p and get
\begin{equation}\tag{7.4}
\frac{1}{2}(\log p)^2+\frac{1}{2}\log p-\frac{(\log p)^2}{2}+\frac{(\log p)^3}{3p}>d\log p+d.
\end{equation}
\begin{equation}\tag{7.5}
d<\frac{\log p}{2(\log p+1)}+\frac{(\log p)^3}{3p(\log p+1)}<\frac{1}{2}-\frac{1}{2\log p}+\frac{1}{\log p(\log p+1)}.
\end{equation}
\end{proof}

\noindent {\bfseries Lemma 3. (Mertens’ third theorem) }
\textit{ For any integer $n>7\,713\,133\,853$, we have
\begin{equation}\tag{L3.1}
\sum_{p\leq n}\log \left(\frac{p}{p-1}\right)=\log \log n+\gamma+R(n),
\end{equation}
where $\gamma$ is the Euler-Mascheroni constant, R(n) is the remainder such that
\begin{equation}\tag{L3.2}
 -\frac{0.005586}{(\log n)^2}<R(n)<\frac{0.005586}{(\log n)^2}.
\end{equation}
}
\begin{proof}
By setting $k=2, \eta_2=0.01$ in Theorem 5.9 of [Dusart 2018], we have, for $n>7\,713\,133\,853$,
\begin{align*}
\lvert R(n)\rvert &<\frac{0.01}{2(\log n)^2}+\frac{4}{3}\cdot \frac{0.01}{(\log n)^3}\\
&=\frac{0.01}{(\log n)^2}\left(\frac{1}{2}+\frac{4}{3\log n}\right)<\frac{0.005586}{(\log n)^2}.\tag{L3.3}
\end{align*}
\end{proof}

\noindent {\bfseries Lemma 4. }
\textit{ Let $g(x)=(\log x)e^{R(x)}$, where
\[R(x)=\frac{0.005586}{(\log x)^2},\]
then $g(x)$ is strictly increasing in interval $(1.1115,\infty)$.}
\begin{proof} Take derivative, we get 
\begin{align*}
 g'(x)&=\frac{1}{x}e^{R(x)}+(\log x)e^{R(x)}\left(-\frac{2\times 0.005586}{x(\log x)^3}\right) \\
&=\frac{1}{x(\log x)^2}e^{R(x)}\left((\log x)^2-0.011172\right).
\end{align*}
So, $g'(x)$ has a zero at $x=1.1115$, and is positive on the right. 
\end{proof}

\noindent {\bfseries Theorem 8. }
\textit{let $N>10^{(10^{13})}$ be an LG number, then
\begin{equation}\tag{8.1}
G(N)<e^\gamma+\frac{0.00995}{(\log \log N)^2}.
\end{equation}
}

\begin{proof} It is easy to see
\begin{equation}\tag{8.2}
\rho (N)=\prod_{i=1}^r\frac{p_i-p_i^{-a_i}}{p_i-1}.
\end{equation}
Because a part is smaller than total, we have 
\begin{equation}\tag{8.3}
\rho (N)<\prod_{i=1}^r\frac{p_i}{p_i-1}\leq \prod_{p\leq p_r}\frac{p}{p-1}
\end{equation}
Substitute n by $p_r$ in (L3.1) of Lemma 3, we get 
\begin{equation}\tag{8.4}
\sum_{p\leq p_r}\log \left(\frac{p}{p-1}\right)=\log \log p_r+\gamma +R(p_r)
\end{equation}
here $R(p_r)$ is the remainder. Take exponential of (8.4), 
\begin{equation}\tag{8.5}
\prod_{p\leq p_r}\left(\frac{p}{p-1}\right)=e^\gamma \log (p_r)e^{R(p_r)}
\end{equation}
We get by (8.3)
\begin{equation}\tag{8.6}
\rho (N)<\prod_{p\leq p_r}\frac{p}{p-1}=e^\gamma \log (p_r)e^{R(p_r)}
\end{equation}
By Lemma 4, $\log (p_r)e^{R(p_r)}$  is increasing, and by Theorem 2, $p_r< \log N$, we can replace $p_r$ with $\log N$.
\begin{equation}\tag{8.7}
G(N)=\frac{\rho (N)}{\log \log N} <\frac{e^\gamma \log (p_r)e^{R(p_r)}}{\log \log N} \leq e^\gamma e^{R(\log N)}
\end{equation}
By Lemma 3,
\begin{align*}
\exp (R(\log N))&<\exp\left(\frac{0.005586}{(\log \log N)^2}\right)\\
&=1+\sum_{k=1}^\infty \frac{(0.005586)^k}{k!(\log \log N)^{2k}}<1+\frac{0.005587}{(\log \log N)^2}.\tag{8.8}
\end{align*} 
So
\begin{equation}\tag{8.9}
G(N)<e^\gamma \left(1+\frac{0.005587}{(\log \log N)^2}\right) <e^\gamma+\frac{0.00995}{(\log \log N)^2}.
\end{equation}
\end{proof}

\noindent {\bfseries Theorem 9. }
\textit{let $N>10^{\left(10^{13}\right)}$ be an LG number. Then\\
1)
\begin{equation}\tag{9.1}
p_r> (\log N)\left(1-\frac{0.005587}{\log \log N}\right).
\end{equation}
Conversely, 2)
\begin{equation}\tag{9.2}
\log N \leq p_r\left(1-\frac{0.005589}{\log p_r}\right).
\end{equation}
}

\begin{proof} Proof by contradiction. Assume $p_r\leq \log N\left(1-\frac{0.005587}{\log \log N}\right)$. It is easy to see
\begin{equation}\tag{9.3}
\rho (N)=\prod_{i=1}^r\frac{p_i-p_i^{-a_i}}{p_i-1}.
\end{equation}
Because a part is smaller than total, we have 
\begin{equation}\tag{9.4}
\rho (N)<\prod_{i=1}^r\frac{p_i}{p_i-1}\leq \prod_{p\leq p_r}\frac{p}{p-1}
\end{equation}
Substitute n by $p_r$ in (L3.1) of Lemma 3, we get 
\begin{equation}\tag{9.5}
\sum_{p\leq p_r}\log \left(\frac{p}{p-1}\right)=\log \log p_r+\gamma +R(p_r)
\end{equation}
here $R(p_r)$ is the remainder. Take exponential of (9.5), 
\begin{equation}\tag{9.6}
\prod_{p\leq p_r}\left(\frac{p}{p-1}\right)=e^\gamma \log (p_r)e^{R(p_r)}
\end{equation}
We get by (9.4)
\begin{equation}\tag{9.7}
\rho (N)<\prod_{p\leq p_r}\frac{p}{p-1}=e^\gamma \log (p_r)e^{R(p_r)}
\end{equation}
By Lemma 4, $\log (p_r)e^{R(p_r)}$  is increasing and by assumption, $p_r\leq C\log N$, where $C:=1-0.005587/\log \log N$, we can replace $p_r$ with $C \log N$.
\begin{equation}\tag{9.8}
\rho (N)<e^\gamma \log (p_r)e^{R(p_r)}\leq e^\gamma \log (C \log N))e^{R(C\,\log N)}
\end{equation}
To get a contradiction, we need to prove 
\begin{equation}\tag{9.9}
e^\gamma \log (C \log N))e^{R(C \log N)}<e^\gamma \log \log N .
\end{equation} 
Cancel $e^\gamma$ and substitute $M:=\log N$, the inequality looks simpler:
\begin{equation}\tag{9.10}
 \log (CM)e^{R(CM)}<\log M.
\end{equation} 
It suffices to prove 
\begin{equation}\tag{9.11}
f(M):=\log (CM)e^{R(CM)}-\log M<0.
\end{equation}
By Lemma 3, \[ R(CM)<\frac{0.005586}{(\log (CM))^2}.\]
Expand the exponential and substituting,
\begin{align*}
f(M)&=\log (CM)\left(\sum_{k=0}^\infty \frac{1}{k!}R(CM)^k\right)-\log M\\
&=\log (CM)\left(1+\sum_{k=1}^\infty \frac{(0.005586)^k}{k!(\log (CM))^{2k}}\right)-\log M\\
&=\log C+\log M+\sum_{k=1}^\infty \frac{(0.005586)^k}{k!(\log (CM))^{2k-1}}-\log M\\
&=\log \left(1-\frac{0.005587}{\log M}\right) +\sum_{k=1}^\infty \frac{(0.005586)^k}{k!(\log (CM))^{2k-1}}\\
&=-\sum_{k=1}^\infty \frac{1}{k}\left(\frac{0.005587}{\log M}\right)^k +\sum_{k=1}^\infty \frac{(0.005586)^k}{k!(\log (CM))^{2k-1}}\\
&=\sum_{k=1}^\infty \left(-\frac{(0.005587)^k}{k(\log M)^k}+\frac{(0.005586)^k}{k!(\log (CM))^{2k-1}}\right).\tag{9.12}
\end{align*}
The summands for $k\geq 2$ are obviously negative. For $k=1$, we have
\begin{equation}\tag{9.13}
-\frac{0.005587}{\log M}+\frac{0.005586}{\log (CM)}=\frac{-0.005587\log C-0.000001\log M}{(\log M)\log (CM)}.
\end{equation}
The difference in numerator decreases when M increases, so we need only to test at $M=(\log 10)\times 10^{13}$, and the difference is $-0.00003<0$. This proves $f(M)<0$ and hence N satisfies (RI) by (9.8), which contradicts to N being LG. 

2) Proof by contradiction. Assume $\log N\leq p_r\left(1+\frac{0.005589}{\log p_r}\right)$. Subsititute (9.2) in to the right side of (9.1), we get
\begin{align*}
\log N \left(1-\frac{0.005587}{\log \log N}\right)&>p_r\left(1+\frac{0.005589}{\log p_r}\right)\left(1-\frac{0.005587}{\log p_r\left(1+\frac{0.005589}{\log p_r}\right)}\right)\\
&>p_r\left(1+\frac{0.005589}{\log p_r}\right)\left(1-\frac{0.005587}{\log p_r}\right)>p_r, \tag{9.14}
\end{align*}
when $p_r>2.3\times 10^{13}$. Hence, N satisfies (RI) by proof of 1). This contradicts to N being LG. 
\end{proof}

\noindent {\bfseries Theorem 10. }
\textit{ Let $N>10^{(10^{13})}$ be an integer, p be the largest prime factor of N. Assume p is the largest prime below $\log N$. If  
\begin{equation}\tag{10.1}
\log N>p_r+\frac{1}{2}\log p+\frac{1}{2}-\frac{1}{2\log p}+\frac{1}{(\log p)(\log p+1)},
\end{equation}
then $G(N)>G(N/p)$.
}
\begin{proof} 
Write $\log N=p+\frac{1}{2}\log p+d$, where d is a to-be determine expression. By Lemma 1, $G(N)>G(N/p)$ if and only if
\begin{align*}
&\frac{1}{2}(\log p)^2+d \log p+d+\sum_{k=1}^\infty \frac{(-1)^{k-1}\left(\frac{1}{2}\log p+d\right)^{k+1}}{k(k+1)p^k}\\
>\frac{1}{2} \log p&+\frac{p(\log p)^2}{2\log N}+\frac{(\log p)^2}{2\log N}+(p+1)\sum_{k=1}^\infty \frac{(\log p)^{k+2}}{(k+2)(\log N)^{k+1}}\tag{10.2}
\end{align*}
Since p is the largest prime below $\log N$, by Proposition 5.4 of {[}Dusart 2018{]}, for $p\geq 89\,693$ we have
\begin{equation}\tag{10.3}
p<\log N\left(1+\frac{1}{(\log p)^3}\right),
\end{equation}
\begin{equation}\tag{10.4}
\frac{p}{2\log N}<\frac{1}{2}\left(1+\frac{1}{(\log p)^3}\right).
\end{equation}
Since $N\geq 10^{(10^{13})}$, $\log N>(\log 10 )\times 10^{13}$, the last terms on both sides of (10.2) are in order of $10^{-13}(\log p)^2$ and can be absorbed by rounding: the numerator 1 in (10.3) was rounded from 0.998. We can concentrate on main terms. For $G(N)>G(N/p)$ it suffices to have
\begin{equation}\tag{10.5}
\frac{1}{2}(\log p)^2+d\log p+d>\frac{\log p}{2}+\frac{(\log p)^2}{2}\left(1+\frac{1}{(\log p)^3}\right).
\end{equation}
Hence
\[d(\log p+1)>\frac{1}{2}\log p+\frac{1}{2\log p}.\]
\begin{equation}\tag{10.6}
d>\frac{\log p+\frac{1}{\log p}}{2(\log p+1)}=\frac{1+\frac{1}{(\log p)^2}}{2\left(1+\frac{1}{\log p}\right)}=\frac{1}{2}-\frac{1}{2\log p}+\frac{1}{(\log p)(\log p+1)}.
\end{equation}
\end{proof}

\noindent {\bfseries Theorem 11. }
\textit{ Let  $N>5040$ be an integer, $p>\log N$ be a prime and p is not a factor of N. If
\begin{equation}\tag{11.1}
\log N<p-\frac{1}{2}\log p+\frac{1}{2}-\frac{1}{\log p+1},
\end{equation}}
then $G(N)>G(Np)$.
\begin{proof} We divide the proof in to two cases.\\
Case 1. $p\geq N$.
\begin{align*}
\frac{G(N)}{G(Np)}&=\frac{\rho(N)\log \log (Np)}{\rho(Np)\log \log N}\\
&=\frac{\log (\log N+\log p)}{\log \log N}\left(\frac{p}{1+p}\right)\\
&\geq \frac{\log (2\log N)}{\log \log N}\left(1-\frac{1}{1+p}\right)\\
&=\left(1+\frac{\log 2}{\log \log N}\right)\left(1-\frac{1}{1+p}\right)\\
&=1+\frac{\log 2}{\log \log N}-\frac{1}{1+p}-\frac{\log 2}{(1+p)\log \log N}\\
&=1+\frac{p\log 2-\log \log N}{(1+p)\log \log N}>1.\tag{11.2}
\end{align*}
Case 2. $p<N$.\\
Write $\log N=p-\frac{1}{2}\log p+d$, where d is a to-be-determined expression. By Lemma 2, $G(N)>G(Np)$ if and only if 
\begin{align*}
&p \sum_{k=1}^\infty \frac{(-1)^k(\log p)^{k+1}}{(k+1)(\log N)^k}\\\
&>-\frac{1}{2}(\log p)^2+d\log p-\frac{1}{2}\log p+d+\sum_{k=1}^\infty \frac{\left(\frac{1}{2}\log p-d\right)^{k+1}}{k(k+1)p^k}.\tag{11.3}
\end{align*}
That is, if and only if
\begin{equation}\tag{11.4}
d(\log p+1)<\frac{(\log p)^2}{1}+\frac{\log p}{2}
-\sum_{k=1}^\infty\frac{(\frac{1}{2}\log p-d)^{k+1}}{k(k+1)p^k}+p\sum_{k=1}^\infty\frac{(-1)^k(\log p)^{k+1}}{(k+1)(\log N)^k}.
\end{equation}
Since $p>\log N$, we can replace $\log N$ with p for all terms with $k\geq 2$.
\begin{align*}
&\quad -\sum_{k=1}^\infty\frac{(\frac{1}{2}\log p-d)^{k+1}}{k(k+1)p^k}+\sum_{k=1}^\infty\frac{(-1)^kp(\log p)^{k+1}}{(k+1)(\log N)^k}\\
&>-\frac{p(\log p)^2}{2\log N}-\sum_{k=1}^\infty\frac{(\frac{1}{2}\log p-d)^{k+1}}{k(k+1)p^k}+\sum_{k=1}^\infty\frac{(-1)^{k-1}(\log p)^{k+2}}{(k+2)p^k}\\
&=-\frac{p(\log p)^2}{2(p-\frac{1}{2}\log p+d)}+\sum_{k=1}^\infty\left(\frac{(\frac{1}{2}\log p-d)^{k+1}}{k(k+1)p^k}+\frac{(-1)^{k-1}(\log p)^{k+2}}{(k+2)p^k}\right).\tag{11.5}
\end{align*}
Consider the sum for $k=2j-1$ and $k=2j$,
\begin{align*}
&-\frac{(\frac{1}{2}\log p-d)^{2j}}{(2j-1)(2j)p^{2j-1}}+\frac{(\log p)^{2j+1}}{(2j+1)p^{2j-1}}-\frac{\left(\frac{1}{2}\log p-d\right)^{2j+1}}{(2j)(2j+1p^{2j}}-\frac{(\log p)^{2j+2}}{(2j+2)p^{2j}}\\
&>\frac{(\log p)^{2j}}{p^{2j-1}}\left(-\frac{\left(\frac{1}{2}-\frac{d}{\log p}\right)^{2j}}{(2j-1)(2j)}+\frac{\log p}{2j+1}-\frac{\left(\frac{1}{2}-\frac{d}{\log p}\right)^{2j}\log p}{(4j)(2j+1)p}-\frac{(\log p)^2}{(2j+2)p}\right)\\
&>\frac{(\log p)^{2j}}{p^{2j-1}}\left(-\frac{1}{2^{2j}(4j^2-2j)}+\frac{\log p}{2j+1}-\frac{\log p}{2^{2j}(8j^2+4j)p}-\frac{(\log p)^2}{(2j+2)p}\right)\\
&>\frac{(\log p)^{2j}}{p^{2j-1}}\left(-\frac{1}{8j^2}+\frac{\log p}{2j+1}-\frac{\log p}{32j^2p}-\frac{(\log p)^2}{(2j+2)p}\right)>0.\tag{11.6}
\end{align*}
So for $G(N)>G(Np)$, it suffices to have
\begin{equation}\tag{11.7}
d\log p+d<\frac{1}{2}(\log p)^2+\frac{1}{2}\log p-\frac{p(\log p)^2}{2(p-\frac{1}{2}\log p+d)}.
\end{equation}
Since
\begin{equation}\tag{11.8}
1-\frac{p}{p-\frac{1}{2}\log p+d}=\frac{-\log p+2d}{2p-\log p+2d}>-\frac{\log p}{2p}>-\frac{1}{(\log p)^2},
\end{equation}
it suffices to have
\begin{equation}\tag{11.9}
d\log p+d<\frac{1}{2}\log p-\frac{1}{2}.
\end{equation}
That is
\begin{equation}\tag{11.10}
d<\frac{\log p-1}{2(\log p+1)}=\frac{1}{2}-\frac{1}{\log p+1}.
\end{equation}
\end{proof}

\begin{center}
{\bfseries \large References}
\end{center}

\noindent {[}Briggs 2006{]} K. Briggs. \textit{Abundant numbers and the Riemann hypothesis}. Experiment. Math., 15(2):251–256, 2006.\\
{[}Broughan 2017{]} K. Broughan, \textit{Equivalents of the Riemann Hypothesis} Vol 1. Cambridge Univ. Press. (2017)\\
{[}CLMS 2007{]} Y. -J. Choie, N. Lichiardopol, P. Moree, and P. Sol\'{e}. \textit{On Robin’s criterion for the Riemann hypothesis}. J. Th\'{e}or. Nombres Bordeaux, 19(2):357–372, 2007.\\
{[}CNS 2012{]} G. Caveney, J.-L. Nicolas, and J. Sondow, \textit{On SA, CA, and GA numbers}, Ramanujan J. 29 (2012), 359–384.\\
{[}Dusart 1998{]} P. Dusart \textit{Sharper bounds for $\psi$, $\theta$, $\pi$, $p_k$}, Rapport de recherche $n\,1998-06$, Laboratoire d’Arithm\'{e}tique de Calcul formel et d’Optimisation\\
{[}Dusart 2018{]} P. Dusart. \textit{Explicit estimates of some functions over primes}. Ramanujan J., 45(1):227–251, 2018.\\
{[}Morrill;Platt 2018{]} T. Morrill, D. Platt. \textit{Robin’s inequality for 25-free integers and obstacles to analytic improvement} \\
https://arxiv.org/abs/1809.10813\\
{[}NY 2014{]} S. Nazardonyavi and S. Yakubovich. \textit{Extremely Abundant Numbers and the Riemann Hypothesis}. Journal of Integer Sequences, Vol. 17 (2014), Article 14.2.8\\
{[}Robin 1984{]} G. Robin. \textit{Grandes valeurs de la fonction somme des diviseurs et hypoth\'{e}se de Riemann}. Journal de mathématiques pures et appliqu\'{e}es. (9), 63(2):187–213, 1984.

\end{document}